# Operational Control of a Multi-energy District Heating System: Comparison of Model-Predictive Control and Rule-Based Control

*M.N.Descamps, N.Lamaison, M.Vallée and R.Bavière*

[a] *Institution, City, Country, e-mail*

**Abstract:**

This study focuses on operational control strategies for a multi-energy District Heating Network (DHN). Two control strategies are investigated and compared: (i) a reactive rule-based control (RBC) and (ii) a model predictive control (MPC). For the purpose of the study a small scale district heating network is modelled using Modelica. The production plant combines a heat pump, a gas boiler and a thermal solar field on the production side with a storage tank for flexibility purposes. On the consumption side, the virtual buildings are aggregated into a single consumer.

We use our co-simulation and control platform, called Pegase, to implement the studied strategies. For both strategies the goal is to meet the consumers' demand while satisfying technical constraints. In addition MPC has the objective to minimize the operational costs, taking into account variable electricity prices and availability of solar thermal resource. Different scenarios are also defined and compared to study the effect of the heat plant sizing and forecasting error.

The operational cost is reduced when switching from RBC to a MPC. As can be expected, MPC is more efficient when dealing with variable energy costs, intermittent solar energy and storage capabilities. This study also demonstrates how our tools enable an easy coupling of Modelica-based simulation with various control strategies. It especially supports the implementation and validation of complex MPC strategies in an efficient way, and yearly simulations are performed within 20 minutes.

# 1 Introduction

Hybrid Energy Networks (HEN), coupling gas, heat and electricity require advanced control strategies to manage the energy generation, storage and consumption, constrained by a high share of intermittent renewable energy sources, and whilst ensuring cost efficiency.

The evaluation and comparison of one particular strategy over an alternative relies on the ability to generate a digital twin of the network considered, and simulate its operation under the chosen strategy.

The digital twin allows generating reproducible and controllable conditions for the operation of the network. This is especially valid when implementing Model Predictive Control (MPC) strategies, which need to be assessed with respect to non-optimal strategies.

The purpose of the current study is to compare the performance of different control strategies applied to a District Heating Network (DHN). The scale and composition of the network are inspired by the existing small scale experimental facility of INES. The benefits of using this simple network are two-fold. First, it is suitable to assess key concepts of district heating systems such as the integration of renewable energy sources, the coupling between the heating network and electricity grid, or the control strategy to reach an optimum operation of the network. Second, this network will be used in the future to validate experimentally these different control strategies. The focus is on the management of the energy production, with flexibility through storage at district level, as opposed to distribution or demand side management, which is not considered here.

The DHN model is operated under two real time control strategies: (i) a reactive rule-based control (RBC) strategy and (ii) a MPC strategy. For each tested strategy, the building heat demand needs to be satisfied at any given instant by the heat production units and the thermal storage discharge.

The paper presents the methodology that is followed to build a DHN digital twin, which relies on simulated and experimental inputs, and operate it on the co-simulation platform PEGASE over a given period. Section 2 describes the case study and the origin of the different components and inputs of the digital twin. The numerical architecture and workflow of the PEGASE platform is briefly introduced in Section 3. The control strategies are presented in Section 4 and applied to the digital twin over a period of approximately 3 months for different scenarios of heat plant sizing. The results are presented and analysed in Section 5.

## 2. Case study

### 2.1. District-Heating network configuration

The heat plant of INES combines a gas boiler and a heat pump, and a thermal solar field produces intermittent heat from RES. All heat production units run in parallel and for this study, we only consider a single consumer. The scale of the DHN is summarised in Table 1.

*Table 1: Assets available on the experimental DHN of INES*

| Component | Capacity |
|---|---|
| Gas boiler | 280 kW |
| Heat pump (condenser side) | 50 kW |
| Consumer peak load | 140 kW |
| Thermal storage | 40 m3 |
| Solar field surface | 140 m2 |

The digital twin assembles the components listed in Table 1 with a similar order of magnitude for the capacities. The core of the digital twin is the distribution network and heat plant simulator which is done with Modelica (Figure 1). Historical measurements coming from the experimental solar field of INES are provided as input to the distribution network and heat plant simulator. The consumer is simulated using TRNSYS.

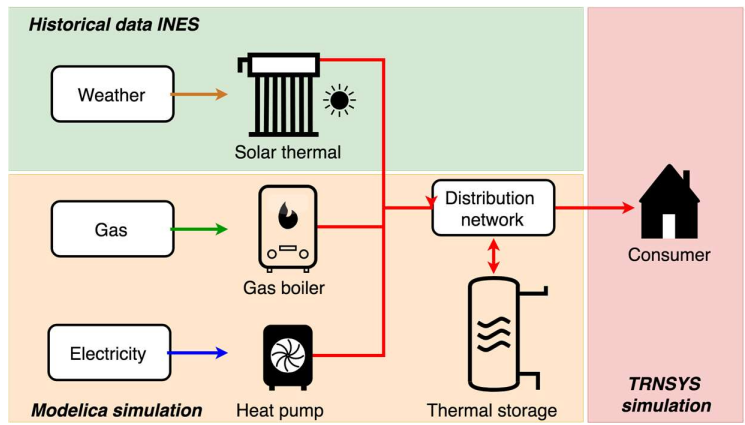
*Figure 1: Heating network sketch and origin of the digital twin components*

## 2.2. Distribution network and heat plant simulator

The heating network is modelled using the equation-based language Modelica along with the simulation platform Dymola. The components of the network are taken from the library Buildings [1] and DistrictHeating [2].

On the production side, the thermal solar field, gas boiler and heat pump are represented each by a simplified model which produces heat from a prescribed value. The input heat values are provided to the model as time series tables. The cold source circuit of the heat pump is not modelled in this work, instead the cold source is assumed constant.

The water tank is represented by a vessel component with perfect temperature stratification. From the heat production side, the hot water enters the tank from the top, and leaves colder from the bottom. On the consumer side, the hot water leaves the top of the tank, and returns colder at the bottom of the tank.

As for the production units, the consumer is represented by a simplified model which consumes heat from a prescribed time-series table, with a negative sign. The pipes are modelled using a node pipe model.

It should be noted that the thermal solar field is curtailed in order to prevent excessive heat supply into the thermal storage in the summer months, when the solar production is high and the consumer demand is low. Experimentally, this corresponds to diverting the solar production away from the thermal storage using a set of valves and an aerothermal refrigerant to dissipate the heat in the environment. Numerically, the curtailment is performed by forcing the solar production to zero when a temperature threshold inside the storage tank is reached.

## 2.3 Thermal solar field

Experimental data was collected for the year 2017 and 2018. The data contains weather information (solar irradiation at 30 degree angle, ambient temperature) as well the output power of a solar field of 140 m2 (Figure 2). The dataset contains outliers and missing values, therefore the current study focuses on a period of continuous and consistent data, which runs from Oct 1$^{st}$ 2017 to Dec 25$^{th}$ 2017.

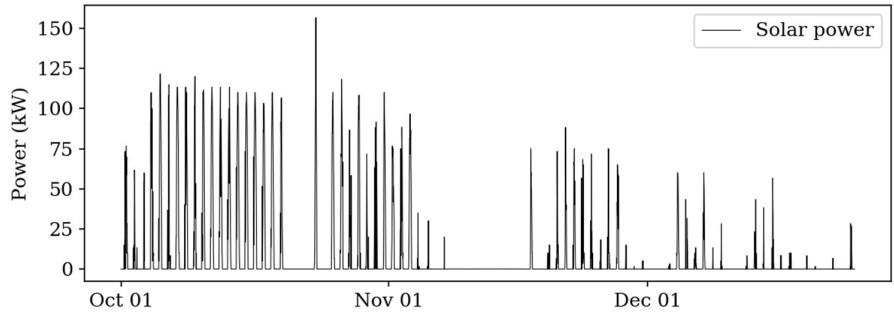

*Figure 2: Experimental measurements of the solar production of INES DHN, Chambery, 2017.*

It is beyond the scope of this paper to develop an accurate prediction model, and interested readers are referred to [3] for the description of machine learning algorithms using data from INES DHN. Instead, a simple prediction correlation is defined by curve fitting the solar production to the solar irradiation on a tilted surface and the ambient temperature. Figure 3 shows the comparison between the predicted and the measured solar production over a period which does not overlap with the curve fit training period.

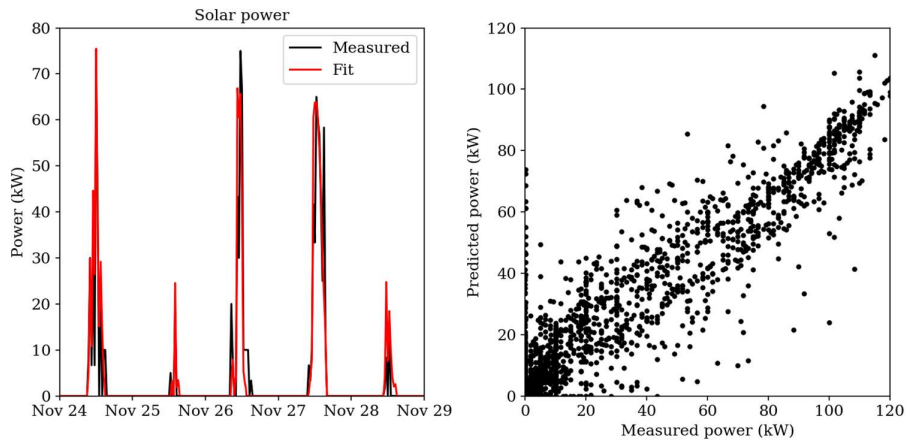

*Figure 3: Prediction of the solar production and comparison to the experimental measurements. Left : period from Nov 24$^{th}$ to Nov 29$^{th}$, Right : period from Oct 1$^{st}$ to Dec 25$^{th}$ 2017*

## 2.2 Consumer heat load

The methodology developed in [4] and [5] is used to simulate the heat load of a reference building with TRNSYS transient system simulation software. This heat load corresponds to space heating and domestic hot water needs. The reference building is chosen as a two story single family building of 140 m² of heated surface.

The heat load profile of the reference building is generated prior to the PEGASE implementation. The simulation is run with a weather file input corresponding to the location of Chambery, France for the year 2017.

The reference building heat load is scaled to match the maximum heat load shown in Table 1. Figure 4 displays the heat load that is applied to the DHN digital twin over the period Oct 1$^{st}$ to Dec 25$^{th}$ 2017.

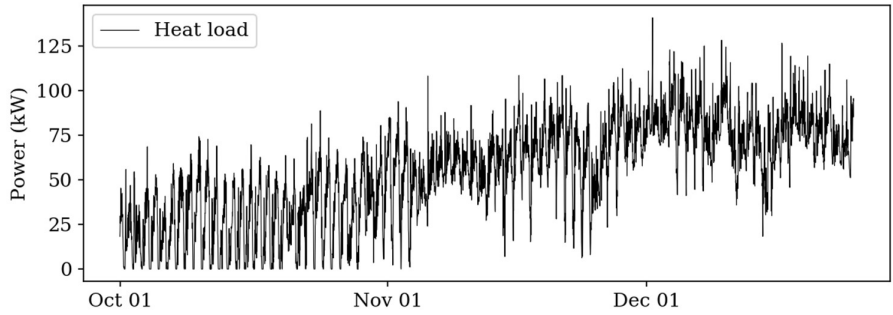

Figure 4: Predefined heat load.

## 2.3 Energy price

The electricity cost corresponds to the year 2017 in France and was obtained using the day ahead prices on the EPEX Spot European Stock exchange published by the association of European Network Operators (ENTSO-E). The gas cost is assumed constant at 0.065 €/kWh following based on fuel cost estimation from [6].. In **Erreur ! Source du renvoi introuvable.** shows the electricity cost per thermal unit assuming a constant coefficient of performance of 3 for the heat pump.

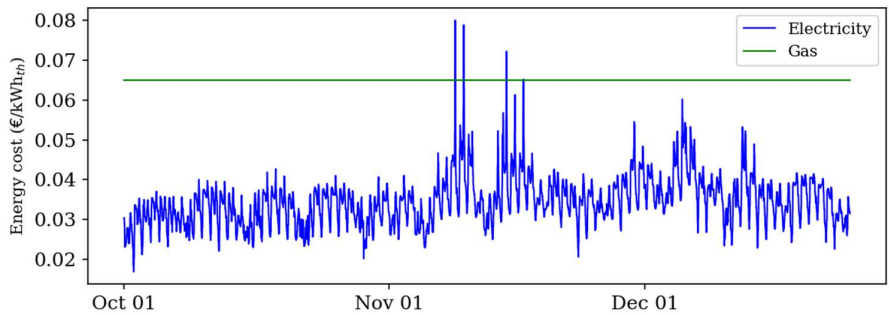

Figure 5: Energy cost per thermal energy unit for the MPC case study

## 2.4. Sizing of the heat plant

The reference scenario (scenario A) considers a heat plant sizing consistent with that of the experimental heating network of CEA-INES. Some variations of scenario A are defined to analyse the effect of heat plant sizing on the control strategies. Scenario B looks at a higher share of heat pump power and lower share of gas boiler power, whereas scenario C looks at a lower share of the solar production by reducing the solar field surface (Table.. ).

Table 2: Heat plant sizing scenarios

| Component | Scenario A | Scenario B | Scenario C |
|---|---|---|---|
| P_GB_max (kW) | 200 | 180 | 200 |
| P_HP_max (kW) | 50 | 70 | 50 |
| Solar field surface (m$^2$) | 70 | 70 | 35 |

# 3 Workflow

The distribution network and heat plant simulator is packaged into a FMU 2.0 (Function Mockup Interface) and run using the PEGASE co-simulation and optimal control environment, alongside two other modules (Figure 6):

- A control module which implements the control strategy RBC or MPC
- In the case of the MPC, a forecast module is required. This module is disabled when running the RBC.

The heating network FMU requires the building heat load and solar production as time series input. After initialisation, an operational cycle runs as follows :

- The control module retrieves the inputs relevant to the chosen strategy and calculates the thermal power of the heat generators. The forecast module is included in this phase if MPC is chosen.
- The calculated heat generator powers are provided as input to the heating network FMU.
- PEGASE retrieves the previous state of the simulation (using the FMU restore state functionality) and launches the heating network simulation for 15 min of physical time
- At the end of the simulation, the FMU saves the internal state of the simulation and provides output for post processing.

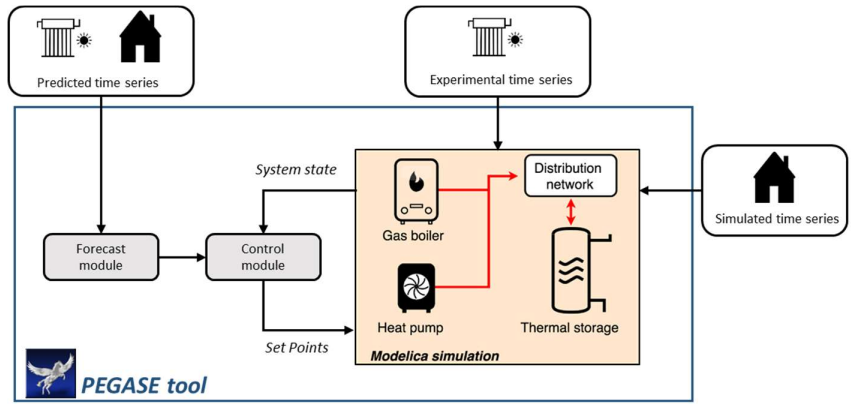

*Figure 6: Architecture of heating network control using the PEGASE environment*

As can be seen, this flexible architecture is suitable for parametric analysis, as it is straightforward to change the modules, for instance replacing RBC by MPC. It is also possible to include hardware components, for instance by replacing the heating network simulation by a real network connected through an Energy Management System.

MPC simulations are performed for perfect weather forecast, i.e. the uncertainty due to weather forecasting error is not taken into account, and the predicted time series are the same as the actual time series except for the solar production. A summary of the time series is found in Table 3.

*Table 3 : List of predefined time series for MPC. Period : Oct $1^{st}$ to Dec $25^{th}$ 2017*

| Time series | Prediction | Actual value |
|---|---|---|
| Weather | Actual weather | Idem |
| Heat load | TRNSYS simulation with weather | Idem |
| Solar production | Curve fit on actual weather file | Experimental measurements |
| Electricity cost | EPEX data for 2017 | Idem |

# 4. Control strategies

The different control strategies are implemented at high level, i.e. by applying power actuators and managing energy constraints. The low level operation of the heating network is performed by the heating network simulator, modifying flow rate and temperature set points to comply with the prescribed powers. As a result, the control strategies are only determining power levels of the different equipment and do not determine temperature levels inside the heating network.

## 4.1 Rule Based Control

The overall control strategy combines the requirement to satisfy the consumer demand to the constraints of maintaining the energy of the thermal storage tank above a minimum value.

The power coming from the solar panels is not a control variable, and is therefore aggregated to the heat consumption (passive components from the operator point of view).

The control strategy maintains the energy of the thermal storage $E$ above a certain threshold $E_{min}$ and imposes a power profile $P_{restore}$ in order to restore the energy in the thermal storage when it is below the chosen threshold. $P_{restore}$ is provided by the injection power of the heat pump $P\_HP$ and gas boiler $P\_GB$, depending on their respective maximum power $P\_HP\_max$ and $P\_GB\_max$.

The RBC gives the priority to the heat pump to provide the base load, and complements the load with

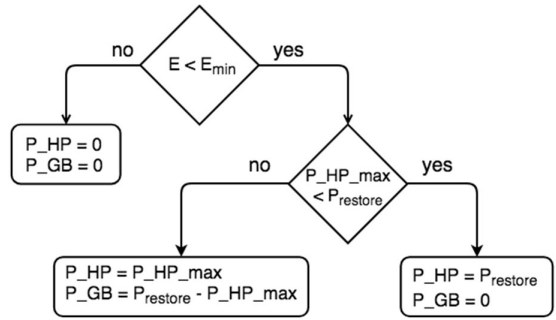

*Figure 7: Rule-based control for the management of the multi-source heat plant and storage*

the gas boiler for peak load, i.e. when the demand is not satisfied only by the heat pump (Figure 7)

The power $P_{restore}$ is proportional to the difference $E - E_{min}$. The proportionality factor depends on the net mass flow rate and a multiplying coefficient chosen to define how fast the energy inside the storage tank should be restored above its threshold. In reality, that proportional coefficient should be linked to the hydraulic capabilities of the installation.

## 4.2 Model Predictive Control

The MPC strategy relies on a Mixed Integer Linear Programming (MILP) formulation to optimize the power flows. The goal is to minimize an objective function while respecting a set of constraints. MILP formulation has been shown to be really adapted to the optimization of energy systems, since it relies on a physico-mathematical approach, is not that costly computationally speaking and represents the optimal solution.

The optimization variables are the power set points of the gas boiler $P\_GB$ and heat pump and $P\_HP$ and the energy of the thermal storage $E$. The optimization problem is the minimization of the cost function which reflects the integral of the operational costs (electricity of the heat pump and gas) over a given time horizon.

In terms of equality constraints, the energy balance between the production and consumption must be respected:

$$\frac{dE}{dt}(t) = P\_HP(t) + P\_GB(t) + Psolar(t) - Pconsumer(t) - K * E(t)$$

where $K$ accounts for heat loss. In terms of inequality constraints, the energy of the thermal storage is bounded, as well as the output power of the production units. Optionally, the first derivatives of the output power can be bounded to prevent too fast unrealistic start-up/shutdown time of the production units.

At each time step, the optimization problem is solved over a calculation horizon of 24h. We use lp_solve as a MILP solver, and if no optimal solution is found at a given time step, the control strategy switches to a default RBC.

## 5. Results

The simulations are performed on laptop computer (16 GB ram; Intel Core i5 2.70 GHz), and runs in around 4 min for MPC and 15 min for RBC for 3 months of simulation time with a time step of 30 min. The longer run time for RBC is due to the event triggering which is inherent to this method and is not computationally efficient.

For RBC and MPC, the share of each heat production unit is shown in Figure 8 over a period of a few days. It can be seen that RBC tends to balance the production and consumption at each time step, with the thermal storage acting as a buffer. Conversely, MPC is not correlated to the instantaneous heat load but is driven by the energy costs, so in Figure 8 the heat pump is switched off during peaks of cost of electricity.

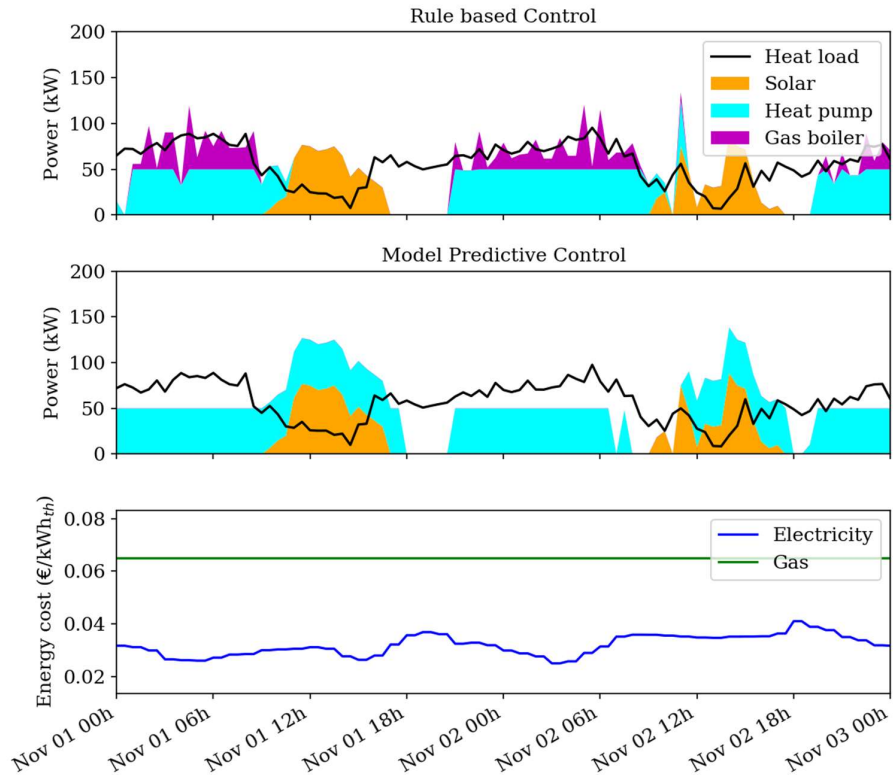

*Figure 8: Share of energy production for scenario A (top and middle) and energy price (bottom)*

The state of charge of the thermal storage has a higher variability for MPC than for RBC. In particular, energy is stored whenever the production costs are low with MPC, and the excess energy is discharged when the heat load is high, as can be seen in Figure 9. In terms of controls, additional constraints could be added to prevent excessive temperatures in the storage tank, which are associated with larger heat loss.

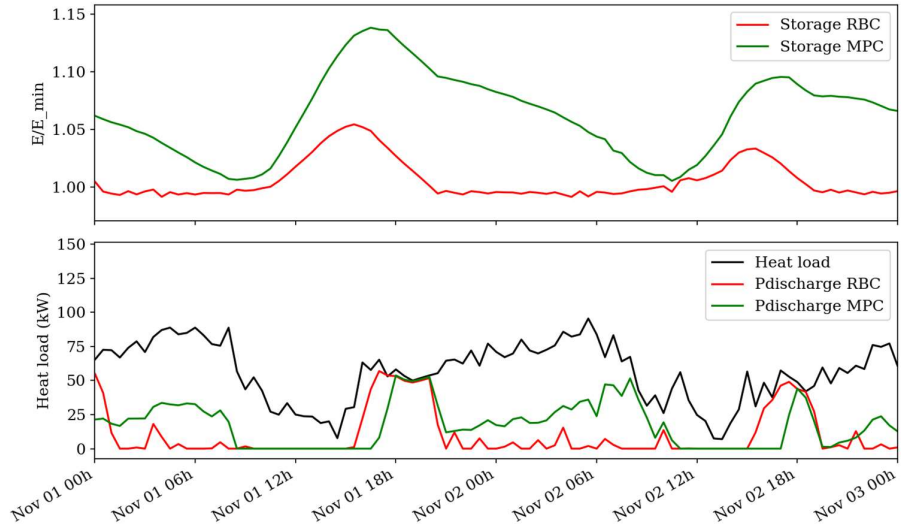

*Figure 9 : State of charge of the thermal storage (top) and associated charge/discharge power (bottom)*

The comparison of both strategies for scenario A is summarized in Table 4 for the period Oct 1st 2017 to Dec 25th 2017. The energy cost is reduced with MPC, which stems from the choice of cost minimization as an objective function. The details for each component show that the main impact of MPC is on the significant reduction of the usage of the gas boiler, to the benefit of the heat pump. This translates into a reduction of energy cost and share of the gas boiler, whereas the same indicators increase for the heat pump. It can be noted that the share of solar energy is slightly lower for MPC than for RBC. A closer inspection of the data shows that this effect is due to prediction errors of the solar power: during sunny days with low heat demand, if the prediction module underestimates the solar power, then excess energy gets stored in the thermal storage. As a result, the maximum temperature of the thermal storage is reached and the distribution network and heat plant simulator curtails the solar power. The performance of the MPC would be improved by adding the solar field into the MILP formulation and including constraint on the temperatures. This would improve the efficiency of the thermal storage and ultimately increase the system flexibility.

*Table 4 Performance of the MPC over RBC for scenario A*

| Indicator | Relative difference MPC-RBC for scenario A |
| --- | --- |
| Total energy cost | -4.6% |
| Energy cost gas | -20.2% |
| Energy cost electricity | +6.5% |
| Energy production | +0.3% |
| Energy share GB | -20.5% |
| Energy share HP | +8.0% |
| Energy share solar | -0.5% |

The influence of changing the heat plant sizing is summarized in Table 5 for both strategies. When increasing the maximum power of the heat pump and reducing that of the gas boiler (scenario B), the heat pump can satisfy the demand for longer periods, which reduces the energy costs for both strategies, because of the lower electricity prices. Similarly, when reducing the solar production (scenario C), the gas boiler needs to be switched on for longer periods, to the detriment of the cost and total energy production. In that case, it is due to the fact the gas boiler gets switched on in night and late evening hours to take advantage of low costs and the energy is stored.

Table 5 Variation of the indicators for scenario B & C compared to scenario A

| Indicator | RBC | | MPC | |
|---|---|---|---|---|
| | Scenario B | Scenario C | Scenario B | Scenario C |
| Energy cost | -13.1% | +18.5% | -14.3% | +20.3% |
| Energy cost gas | -87.7% | +89.1 | -117.0% | +119.5% |
| Energy cost elec | +23.9% | -15.6% | +21.5% | -14.2% |
| Energy production | 0% | 0% | +0.1% | -0.1% |
| Energy share GB | -87.7% | +89.1% | -117.0% | +119.5% |
| Energy share HP | +23.8% | -14.2% | +21.9% | -13.6% |
| Energy share solar | 0% | -66.7% | -0.1% | -66.4% |

## 6. Conclusion and perspective

In this paper, high level control strategies are tested on a digital twin of INES DHN. The methodology relies on the modelling of a DHN in Modelica, which is integrated with other modules (control, forecast) into the co-simulation platform PEGASE.

The system considered is a small scale heating network with a thermal storage tank and a multi-source heat plant consisting of a heat pump, a gas boiler and a solar field. Predefined heat loads are applied to the system, and the control strategies manage the operation of the heat plant over a period of a year.

With the objective to minimize the energy costs, MPC outperforms RBC in all cases. A detailed interpretation of the energy indicators underlines the potential improvements on the share of solar thermal energy. A parametric analysis was carried out to study the influence of the heat plant sizing with respect to the control strategies. Such an analysis illustrates the benefits of using the PEGASE platform, both in terms of run time and flexibility of architecture.

An extension of this work is to replace the simulated DHN by the experimental micro DHN of CEA-INES to evaluate in real time the performance of the control strategy.

## Acknowledgments

The authors gratefully acknowledge the financial support of the funded PENTAGON project by the European Commission, granted in the Horizon 2020 no.731125.

## References

[1] M. Wetter, W. Zuo, T. Nouidui and X. Pang, "Modelica Buildings Library," *Journal of Building Performance and Simulation,* vol. 7, no. 4, pp. 253-270, 2014.

[2] L. Giraud, R. Bavière, M. Vallée and C. Paulus, "Presentation, Validation and Application of the DistrictHeating Modelica Library," in *Proceedings of the 11th International Modelica Conference*, Versailles, France, 2015.

[3] M. Ahmad, J. Reynolds and Y. Rezgui, "Predictive modelling for solar thermal energy systems: A comparison of support vector regression, random forest, extra trees and regression trees," *Journal of Cleaner Production,* vol. 203, pp. 810-821, 2018.

[4] R. Heimrath and M. Haller, "Project Report A2 of Subtask A : the Reference Heating System, the Template Solar System," IEA SHC Task 32, 2007.

[5] U. Jordan and K. Vajen, "DHWcalc: Program to generate domestic hot water profiles with statistical means for user defined conditions," in *Proc. ISES Solar World Congress*, Orlando (US), 2005.

[6] IRENA, "Renewable Energy in District Heating and Cooling : a Sector Roadmap for Remap," 2017.